\documentstyle[12pt]{article}

\def\cp{{\bf CP}^n} \def\g{{\bf g}}
\title{Quantum Sphere via Reflection Equation Algebra}
\author{D.~Gurevich,\\
{\small \it ISTV, Universit\'e de Valenciennes,
59304 Valenciennes, France}\\
P.~Saponov\\
{\small\it Theory Division, Institute for High Energy
Physics, 142284 Protvino, Russia}}

\def\Lq{{\cal L}_q(R)}

\def\Tq{{\cal T}_q(R)}
\begin{document}
\maketitle
\begin{abstract}
Quantum sphere is introduced as a quotient of the so-called Reflection
Equation Algebra. This enables us to construct some line bundles on it
by means of the Cayley-Hamilton identity whose a quantum version was
discovered in \cite{PS}, \cite{GPS}. A new way to introduce some
elements of "braided geometry" on the quantum sphere is discussed.
\end{abstract}


\renewcommand{\theequation}{{\thesection}.{rabic{equation}}}
\def\r#1{\mbox{(}\ref{#1}\mbox{)}}
\setcounter{equation}{0}
\def\uqs{ U_q(sl(2))}
\def\uqsm{ U_q(sl(2))-Mod}
\def\ot2{\otimes 2}
\def\ovalp{\overline{\alpha}}
\def\ovb{\overline{\beta}}
\newcommand{\gggg}{g}
\newcommand\lbc{[}
\newcommand\rbc{]}
\newcommand{\Pm}{P_-(t)}
\newcommand{\rr}{\rho^{\ot 2}}

\def\ot{\otimes}
\def\kEi{k(E_{\la_i})}
\def\Ei{E_{\la_i}}
\def\El{E_{\la}}
\def\Eliq{E_{\la_i}^q}
\def\kMl{k(\Ml)}
\def\oM{\overline{M}}
\def\Mn{M_{\nu}}
\def\oP{\overline P}
\def\RVM{{\cal R}{(V,\, {M_{\mu}^q)}}}
\def\RVVM{{\cal R}{(V^{\ot 2},\, {M_{\mu}^q)}}}
\def\RVMh{{\cal R}{(V,\, {M_{\mu}^{q\h})}}}
\def\kMm{k({\overline M}_{\mu})}
\def\RVm{{\cal R}(V, {\overline M}_{\mu})}
\def\Rn{{\cal R}_{\nu}}
\def\Ru{{\cal R}_{\nu}}
\def\Rmi{{\cal R}_{\mu_i}}
\def\kHq{k(H_q)}
\def\la{\mu}
\def\Ml{{M_{\la}}}
\def\Mlq{{M_{\la}^q}}
\def\co{\cal O}
\def\M{\cal M}
\def\N{{\cal N}}
\def\E{{\cal E}}
\def\gah{{\gamma_h}}
\def\aq{{\cal A}_q}
\def\sq{S^2_q}
\def\ahc{{\cal A}_h^c}
\def\ahq{{\cal A}_{h,\,q}}
\def\usl{U(sl(2)}
\def\ush{U(sl(n)_{\h}}
\def\Uh{U(gl(n)_{\h})}
\def\Uhq{U(gl(n)_{\h,q})}
\def\uslh{U(sl(2)_h}
\def\uq{U_q(sl(2))}
\def\Uq{U_q(sl(n))}
\def\Ush{U(sl(n)_{\h})}
\def\Ushq{U(sl(n)_{\h,q})}
\def\aam{A^{(m)}}
\def\aan{A^{(n)}}
\def\vn{V^{\ot n}}
\def\iin{I^{(n)}}
\def\inn{I^{n}}
\def\de{\delta}
\def\De{\Delta}
\def\RR{\bf R}
\def\CC{\bf C}
\def\cp{{\bf CP}^n}
\def\Mat{\rm Mat}
\def\Sym{\rm Sym}
\def\ahq{A_{\h,\,q}}
\def\om{\Omega}
\def\omo{\Omega^1}
\def\omd{\Omega^2}
\def\aqc{{\cal A}^c_{ 0\, q}}
\def\ahqc{{\cal A}^c_{\h\, q}}
\def\ac{{\cal A}^c_{0\, 1}}
\def\h{{\hbar}}
\def\gh{{\g_{\hbar}}}
\def\ah{{\cal A}_{\hbar}}
\def\V{\bf V}
\def\VV{{\bf V}^{\ot 2}}
\def\vv{{ V}^{\ot 2}}
\def\ks{k(S^2)}
\def\ksq{k(S^2_q)}
\def\thq{T(H_q)}
\def\am{A_{\mu}}
\def\ami{A_{\mu}^i}
\def\gq{{\g}_q}
\def\ggq{{\g}_q^{\ot 2}}
\def\sss{\cal S}
\def\Ob{Ob\,}
\def\h{{\hbar}}
\def\a{\alpha}
\def\na{\nabla}
\def\ve{\varepsilon}
\def\ep{\epsilon}
\def\bea{\begin{eqnarray}}
\def\eea{\end{eqnarray}}
\def\beq{\begin{equation}}          \def\bn{\beq}
\def\eeq{\end{equation}}            \def\ed{\eeq}
\def\nn{\nonumber}
\def\oH{{\overline H}}
\def\A{{\cal A}}
\def\aaa{\cal A}
\def\rank{{\rm rank\, }}
\def\End{{\rm End}\,}
\def\Ker{{\rm Ker\, }}
\def\root{{\rm root\, }}
\def\diag{{\rm diag\, }}
\def\Im{{\rm Im\, }}
\def\Id{{\rm id }}
\def\ra{{\rm rank\, }}
\def\Id{{\rm id\, }}
\def\Vect{{\rm Vect\,}}
\def\Hom{{\rm Hom\,}}
\def\uHom{\underline{\rm Hom}\,}
\def\Tr{{\rm Tr}}
\def\det{{\rm det}}
\def\codet{{\rm codet}}
\def\dim{{\rm dim}\,}
\def\id{{\rm id}\,}
\def\span{{\rm span}\,}
\def\udim{\underline{\rm dim}\,}
\def\dem{{\rm det}^{-1}}
\def\Fun{{\rm Fun\,}}
\def\de{\delta}
\def\al{\alpha}
\renewcommand{\theequation}{{\thesection}.{\arabic{equation}}}
\setcounter{equation}{0}
\newtheorem{proposition}{Proposition}
\newtheorem{conjecture}{Conjecture}
\newtheorem{corollary}{Corollary}
\newtheorem{theorem}{Theorem}
\newtheorem{definition}{Definition}
\newtheorem{remark}{Remark}
\def\R{{\cal R}}
\def\n{{\bf n}}   \def\va{\varepsilon}
\section{Introduction}
Quantum sphere $\sq$ is the simplest example of a quantum variety. It
was introduced
in \cite{P}
by means of a quantum reduction procedure which in a schematic way
can be represented  by the formula
$$\sq=SU_q(2)/(SU_q(1)=SU(1)).$$
This means that the quantum sphere (by abusing the language we use the
term
 "quantum variety" for the corresponding coordinate rings\footnote{All
 coordinate rings are treated
in the spirit of affine algebraic geometry as  quotients of a
polynomial ring. The basic field $k$
is assumed to be $\bf C$ or $\bf R$. In the latter case $q$ is assumed
to be real. The choice of the field is similar to the classical case.
In any case $q$ is
generic.}) is realized as
a subspace in the algebra $k(SU_q(2))$.
Other quantum varieties  are usually introduced in a similar way.
Nowadays, the language of
the so-called Galois-Hopf extension is employed. It makes use of a
couple of Hopf algebras giving rise
to a "quantum coset" (cf. \cite{S}, \cite{BM}).

However, this approach has some defects.
First, it does not allow us to control flatness of deformation of
reduced objects\footnote{We
refer the reader to \cite{DGK} for a definition. However, if a
classical
(resp., quantum) object decomposes into a direct
sum of irreducible  $U(\g)$- (resp., $U_q(\g)$-) modules the
flatness means that the both objects consist
of similar components.}. It is worth noticing that the differential
calculus on
$SL_q(n)$ constructed by S.Woronowicz \cite{W1}, \cite{W2}
with the use of the Leibniz rule is not a flat deformation
of its classical counterpart (cf. \cite{AAM},  \cite{Ar}, \cite{FP1},
\cite{FP2}, \cite{I}).
However, this
differential calculus plays the crucial role in all constructions 
(metric, bundles, connection etc.) of
quantum (braided) version
of differential geometry  on quantum
varieties. So, these
objects are pointless in a sense, if we want to treat them as an
approximation of their classical counterparts.

Second, the above way of introducing quantum varieties
 cannot be generalized to "non-quasiclasical quantum varieties",
 associated to
non-standard quantum $R$-matrices constructed in \cite{G}\footnote{In
the sequel we use the term
quasiclassical for the objects related to the quantum group $\Uq$.}.
Although a $k(G_q)$ type algebra can be apparently defined for
numerous non-standard $R$-matrices (cf. \cite{G}), it does not have
any subalgebra which
would give rise to a "quantum coset".

The main purpose of this note is to present another way of  defining
"quantum varieties"
which on one hand would be
valid both in the quasiclassical and non-quasiclassical case and on
the other hand would enable us
to control flatness of deformation.
In the framework of this approach the quantum function algebras or
their dual objects
(quantum groups) play only the role of symmetry groups (however,
we can do without  these objects at all). The crucial role in our
approach belongs to the
so-called reflection equation (RE) algebra. We consider this algebra
as a very fruitful tool of
"braided geometry".

The central objects of this geometry are quantum (braided) varieties
which are introduced as some
quotients of the RE algebra. Otherwise stated, they are realized
explicitly by means of some "braided
system of equations" in the spirit of affine algebraic geometry. In
virtue of \cite{PS},
\cite{GPS} there exists some polynomial (of Cayley-Hamilton type)
identity
for the matrix formed by the generators of this algebra. Moreover, the
coefficients of the
corresponding polynomial are elements of the center of the RE algebra.
This property enables us to introduce  quantum line bundles on some
quantum varieties
(in particular, quantum sphere) in terms of projective modules in the
spirit of
the Serre-Swan approach (cf. \cite{Se},
\cite{Sw}) and to show the projectivity of the corresponding modules.

In this paper we
want to demonstrate the
usefulness and the  power of the described approach in "braided
geometry"  on the example of quantum sphere. The note is organized
as follows.  In section 2 we introduce quantum sphere (or what is the
same quantum hyperboloid if we
ignore involution operator) and discuss a way to introduce
differential calculus,
tangent space and some other structures on it without either
reduction procedure or  Leibniz
rule at all. Quantum group (or their dual objects) are used only as
substitutes of symmetry groups.
In section 3 we introduce the RE algebra and describe some its
properties.
In  section 4 we realize the quantum sphere in terms of this algebra
and use this
realization in order to introduce line bundles on quantum sphere via
the Cayley-Hamilton identity.

{\bf Acknowledgment} The authors are supported by the grant
PICS-608/RFBR 98-01-22033.

\section{Differential calculus and tangent space on quantum sphere}
The purpose of this section is to describe a way to introduce a
quantum
sphere
(hyperbo\-loid)\footnote{If the coordinate ring of a sphere or a
hyperboloid
is defined as a quotient of the polynomial ring in the spirit of
affine algebraic geometry
no difference between
these varieties appears.
In particular, all irreducible $sl(2)$-modules involved in
differential calculus are
finite-dimensional.
A difference appears while we consider the hyperboloid equipped with
other functional spaces. For the same reason we do not make any
difference between
 the QG $\Uq$  and $U_q(su(2))$.}
and some derived  objects
of quantum geometry  (differential calculus, tangent space)  on it
without any reduction procedure.

Let $V$ be spin 1 $\uq$-module. Then the space $\vv$ being equipped
with a structure of
$\uq$-module decomposes into a direct sum of three irreducible
components
\beq\
\vv=V_0 \oplus V_1 \oplus V_2. \label{fusion}
\eeq
where $V_i$ stands for the spin $i$ $\uq$-module. (We avoid using any
coordinate writing.)

Then a quantum sphere (or more precisely, its the coordinate ring)
$\ksq$ can be defined as the
quotient
\beq
T(V)/\{V_1,\, v_0-c\},\,\, c\in k\,\,{\rm is \,\, assumed \,\,to\,\,
be\,\, fixed}.\label{dom}
\eeq
Here $T(V)$ stands for the free tensor algebra of the space $V$,
$\{I\}$ stands for the ideal
in $T(V)$ generated by a subset $I$ and  $v_0$ is a generator of the
1-dimensional component $V_0$.

We consider this quotient as a "q-commutative" algebra. Its
"q-non-commutative" counterpart can
be defined in a similar way if we replace $V_1$ in the denominator of
\r{dom} by
$V_1-\h\al(V_1)$ where
$\al: \vv\to V$ is a non-trivial $\uq$-morphism (it is unique up to a
factor which in the sequel is
supposed to be fixed). Thus,  we get two parameter algebra denoted
$\ahq$ which becomes
$\ksq$ as $\h=0$ and moreover it is a flat deformation of the
classical
counterpart
$A_{0,1}=k(S^2)$. Let us note that the Podles' sphere is another
parametrization of this algebra
equipped with an involution operator (see a discussion on an
involution operator in section 3).

\begin{remark} For the Lie algebras $\g=sl(n),\,\,n>2$ the following
difficulty appears.
Let us equip $\g$ considered as a vector space with the structure of a
 $\Uq$-module. Then it is no clear
what is a reasonable way to define a morphism analogous to that $\al$
above since
in the space $\g^{\ot 2}$ the component isomorphic to $\g$ itself
occurs twice.
By the same reason it is not clear what is a q-analogue of the
symmetric algebra of the space $\g$.
This problem can be solved in terms of the so-called reflection
equation algebra
(see  section 3).
\end{remark}

\begin{remark}
Let us remark that there exist other tensor or quasitensor categories
whose fusion rings
look like that of category of
$sl(2)$-modules. Let $R:\VV\to\VV$ be a Hecke symmetry, i.e., a
solution of the Yang-Baxter
(YB) equation
satisfying the quadratic equation
$$R^2=\id+(q-q^{-1})R$$
 such that the Poincar\'e series of the corresponding
"skewsymmetric algebra" (cf. \cite{G}) is of the form
$$P_-({\V}, t)=1+nt+t^2,\,\, n=\dim{\V}.$$

A big family of such type Hecke symmetries was constructed in
\cite{G}.

Then in the category generated by the space $\V$ in the spirit of the
paper \cite{GM} there is a
space $V$ whose tensor square decomposes similarly to \r{fusion}.

The non-quasiclassical quantum (or braided) variety looking like
the quasiclassical
quantum sphere
and its non-commutative counterpart
can be constructed in the same way as above without any reduction
procedure.
Also remark that  there exists
another quantum sphere (hyperboloid) being a deformation of the
classical
one. This quantum sphere
corresponds to
the involutary solution of the quantum YB equation arising from the
triangular
classical r-matrix $H\wedge X$.
\end{remark}

In a similar non-coordinate manner and without any reduction procedure
there
can be introduced spaces of differentials  on the quantum sphere.
It is done in \cite{A}, \cite{AG}.
Let us describe briefly a way to introduce
$\omo$ suggested in these papers. Denote $V'$ the space isomorphic to
$V$ as $\uq$-module, but
generated by the differentials of the elements from $V$. Let
$\ksq\ot V'$ be free left $\ksq$-module. Let us introduce
the space $\omo$ as its quotient over the submodule
generated by the elements $(V\ot V')_0$.  Hereafter we use the
notation
$(V_i\ot V_j)_k$ for the spin $k$ component in the tensor product
$V_i\ot V_j$.

The second differential space $\omd$ can be introduced similarly to
$\omo$. It is not
difficult to describe  decompositions of the $\ksq$-modules $\omo$ and
$\omd$ and of the algebra
$\ksq$ itself
into direct sums of $\uq$-modules  for a generic
$q$ (cf. \cite{AG}). They look like the decompositions of
their classical counterparts into sums of $sl(2)$-modules.

Let  $d$ be the differential acting in the classical differential
algebra
(as usual, $d$ is subject to the Leibniz rule).
Since $d$ is a $\uq$-morphism it maps any
$\uq$-module either to an isomorphic module or to 0. Let us define
$d$ in the quantum case in a similar way by sending any $\uq$-module
containing in $\ksq$-module $\omo$, $\omd$ or the algebra $\ksq$
itself  either to an isomorphic $\uq$-module or to 0 similarly to the
classical case (also we can require the quantum operator $d$ to be a
deformation of its classical counterpart).

 By this the quantum differential $d$ is defined uniquely up to a
 factor  on any irreducible
$\uq$-component. By construction we have just the same cohomology as
in the classical case:
$$\dim\, H^0=1,\,\, \dim\, H^1=0,\,\, \dim\, H^2=1.$$
Moreover, the same result will be valid if we realize this
construction on any
non-quasiclas\-sical sphere mentioned in remark 2.

Thus, we have constructed our version of quantum differential calculus
without Leibniz rule and without any algebraic structure in the space
$$\om^*=\Omega^0\oplus\omo\oplus\omd,\,\, \Omega^0=\ksq.$$

Let us pass now to discussion of what the tangent bundle on the
quantum sphere is
(or, equivalently, what the phase space would be if we considered the
quantum sphere as a configuration space).
This problem can be split into two parts. First, we want to represent
the quantum tangent space as
a $\ksq$-module. Second, we want to assign to the elements of this
space an operator meaning.

Let us begin with the classical case. The tangent space on the usual
sphere can be introduced by  the equation
\beq
(V\ot V')_0=0\;, \label{tan}
\eeq
where the space $V'$ is generated not by differentials as it was the
case for the quantum sphere but
by the elements of the Lie algebra $su(2)$ itself
(we consider here the compact form of the Lie algebra in question).
In the explicit coordinate form the above equation looks as follows
\beq
xX+yY+zZ=0 \label{tanspace}
\eeq
where $x,\,y,\,z$ are generators of the algebra $\ks$ and $X, \,Y,\,
Z$  are the corresponding
generators of the Lie algebra $su(2)$. As operators they can be
represented by infinitesimal
rotations:
$$X=z\partial_y-y\partial_z,\,\,
Y=x\partial_z-z\partial_x,\,\,Z=y\partial_x-x\partial_y.$$

On the hyperboloid the relation \r{tan} takes the form
\beq
xY+yX+{1\over 2} hH=0, \label{tann}
\eeq
where $X,\, Y,\,H$ are "hyperbolic infinitesimal rotations".

Similarly, for any regular variety $M$ the tangent space $T(M)$
has a structure of a $k(M)$-module and reciprocally $k(M)$ has a
$T(M)$-module structure
since $T(M)$ acts on $k(M)$.
These two structures are compatible in some sense  giving rise to the
notion
of Lie-Rinehart algebra
(cf. \cite{R}).

What is a q-analogue of the space $T(S^2)$ or otherwise stated what
are q-analogues of
the operators $X,\, Y,\,H$ above which would generate the quantum
tangent space? The operators
$X, \, H,\, Y$ coming from the quantum group $\uq$ do not fit for this
role since
they do not satisfy any relation which would be a deformation of
\r{tann}.
As a $\ksq$-module the quantum tangent space can be introduced by
relation \r{tan} considered in the corresponding tensor or quasitensor
category (the same is valid in the non-quasiclassical case). However,
it is not evident in advance whether there exist operators defined on
the algebra $\ksq$ and satisfying this relation. Nevertheless, as was
shown in \cite{A} there exist some operators satisfying relations
\r{tan} and the defining relations of the algebra $\ahq$. We will
treat these operators and all their linear combinations with
coefficients from $A_{0,q}=k(S^2_q)$ as q-analogues of vector fields
on the quantum  sphere (hyperboloid). By definition, the quantum
tangent space $T(S^2_q)$ being a $k(S^2_q)$-module and equipped with a
$\uq$-covariant action
$$
T(S^2_q)\ot k(S^2_q)\to k(S^2_q)
$$
is just the space of all vector fields.

Suggested in \cite{A} (also cf. \cite{AG}) was a method of
constructing some q-analogues of metric and connection in terms of
these quantum (braided) vector fields similarly to the classical way
of introducing these operations (but without using any form of Leibniz
rule).

\section{Reflection equation algebra}
\setcounter{equation}{0}

The purpose of this section is to introduce the RE algebra and to
describe its properties.

Note that this algebra in the case of a quantum R-matrix depending on
a spectral parameter was introduced by Cherednik as a boundary
condition in the inverse
scattering problem. Later some different versions of this algebra
 were studied in a similar context
in \cite{KSk}, \cite{KSa}. In the framework of braided geometry it was
introduced by Sh.Majid
\cite{M1} (cf. also \cite{M2}). Also he has shown that this algebra
possesses a braided Hopf algebra structure.
(In the case of an involutary operator $R$ a similar Hopf algebra
was introduced earlier by one of the author, cf. \cite{G} and
references therein).

As RE algebra we call the algebra generated by $n^2$ elements
$l_i^j,\,\, 1\leq
i, \, j\leq n$
 subject to the following matrix relation
\beq
RL_1RL_1=L_1RL_1R \label{RE}
\eeq
where
$$
L_1=L\ot \Id\,\, {\rm and}\,\, L=(l_i^j). $$
This algebra will be denoted  $\Lq$.

If $R$ is a  Hecke symmetry and is a deformation of usual flip
then the RE algebra is a flat deformation of its classical
counterpart which is the symmetric algebra of the linear space
$\span(l_i^j)$. This was shown in \cite{L} (or by another method in
\cite{D1}).

Let us describe some other properties of this algebra.
First of all, in the quasiclassical case
its product is $\Uq$-covariant.  Habitually, this property is
presented in a
dual form. Namely, the RE algebra is equipped with a left RTT-comodule
structure
so that
the product in the RE algebra is covariant w.r.t. the coaction
\beq
\de: \Lq\to\Tq\ot \Lq,\,\,\, \de(l_i^j)\to  t^j_p\,S(t_i^k)\ot l_k^p
\label{copr}
\eeq
(a summation upon repeated indices is assumed).
Hereafter $\Tq$ stands for the famous RTT algebra (cf. \cite{FRT})
defined by the relation
$$
RT_1T_2=T_1T_2R\;,
$$
and $S$ is the antipode in it.
(Note that this property is valid in the
non-quasiclassical case if the algebra $\Tq$ is well defined, cf.
\cite{G}.)

Second, the RE algebra being a quadratic algebra admits a
quadratic-linear counterpart which is its flat deformation and looks
like the enveloping algebra $U(gl(n))$. More precisely, this
quadratic-linear deformation (denoted $\Uhq$) tends to the RE algebra
as $\h\to 0$ and to $U(sl(n)_{\h})$ as $q\to 1$. Here we use the
notation $\g_{\h}$ for the  Lie algebra which differs from $\g$ by the
factor $\h$ introduced in the  bracket of the algebra $\g$. By this we
want to represent the algebra $\Uh$ as a deformational object w.r.t.
the symmetric algebra $\Sym(\g)$.

Third, this algebra unlike the RTT algebra, has a big center. In
particular, the so-called quantum trace $\Tr_q(L)$ is central (cf.
\cite{FRT}). By contrast, a similar trace in the RTT algebra is not
central. Moreover, for the RE algebra some quantum analogues of the
Newton relations and the  Cayley-Hamilton (CH) identity hold. It was
shown in \cite{GPS} (in the quasiclassical case these relations were
previously established in \cite{PS}). Furthermore, some version of the
CH identity is also valid for the algebra $\Uhq$ (cf. \cite{GS}). By
passing to the limit $q \to 1$ we get a similar identity for the
algebra $\Uh$.

Let us make a precise.  For the matrix $L$ formed by generators of one
of the above algebras ($\Lq$, $\Uhq$, or   $\Uh$) the following
relation is valid
\beq
\sum^p_{i=0}(-L)^i\sigma(p-i)=0\;.
\label{CH}
\eeq
Here $p$ is the rank of the Hecke symmetry $R$ (cf. \cite{G}, in the
quasiclassical case $p=n$) and the coefficients $\sigma(i)$ are
central elements of the algebra in question.

Following \cite{GS} we introduce "quantum generic orbits" as the
quotients of the RE algebra over the ideal generated by the elements
$$
\{\sigma(i)-c_i,\,\,c_i\in k\}\;.
$$

By changing $\sigma(i)$ for the numbers $c_i$ in relation \r{CH} we
get a polynomial $\oP$ with numerical coefficients. This polynomial
plays a crucial role in defining line bundles over quantum varieties
(see  section 4).

Also observe that the operators
$$
\de_1:l_i^j\to l_k^j\, l_i^k,\,\, \de_2:l_i^j\to
l_m^j\,l_k^m\,l_i^k,...\
$$
which map the space $\span (l_i^j)$ to the RE algebra are compatible
with coaction \r{copr}.
We will call such maps {\em morphisms}.
Otherwise stated, this property means that the maps
$L\to L^k,\,\, k = 2,\,3,...$ are morphisms.

Namely, this property allows us to ensure compatibility of identity
\r{CH} with coaction \r{copr} (or the action of the QG $\Uq$ in the
quasiclassical case). Note that the images of the above operators
belong to the RE algebra itself but not to its tensor powers. These
operators are equal to the product of the comultiplication and
multiplication operators.

Remark that for the RTT algebra the maps $T\to T^k$ do not possess
similar property. For this algebra there exists some analogue of the
CH identity \cite{IOP} but it is not as nice as \r{CH} is. First, the
role of the powers  $T^k$ is played by some complicated enough and
less natural expressions and second, the coefficients of the
corresponding polynomial are not central. It is this fact that
prevents us from introducing quantum orbits as some quotients of RTT
algebra. For this algebra only reduction procedure is admitted.
However, the quantum sphere can be realized in the both ways: as a
quantum coset and as a restriction of the RE algebra.

Let us explain this at the quasiclassical level, i.e. by means  of
Poisson brackets. The quasiclassical counterpart of the RTT algebra is
well known. It is the so-called Sklyanin bracket. It can be reduced to
any semisimple orbit in $sl(n)^*$. However, it is not defined on the
whole $sl(n)^*$. While the Poisson bracket being the quasiclassical
counterpart of the RE algebra (denoted $\{\,\,,\,\,\}_{RE}$) is
well-defined on $gl(n)^*$ (and $sl(n)^*$ if we kill the trace).
Moreover, this bracket is compatible with the linear Poisson-Lie
bracket (denoted $\{\,\,,\,\,\}_{PL}$), i.e. they form the Poisson
pencil
$$
\{\,\,,\,\,\}_{ab}=a\{\,\,,\,\,\}_{PL}+b\{\,\,,\,\,\}_{RE},\,\,\, a,
b\in k.
$$
The simultaneous quantization of this pencil is just the algebra
$U(gl(n)_{\h,q})$ mentioned above.

This Poisson pencil can be restricted on any semisimple orbit (as well
as on any other type of orbit, cf. \cite{D2}) in $gl(n)^*$ (or
$sl(n)^*$). So, on such an orbit we have the reduced Sklyanin bracket
and the restricted Poisson pencil $\{\,\,,\,\,\}_{ab}$. In general,
the reduced Sklyanin bracket has nothing in common with this pencil
but on a symmetric orbit is becomes a particular case of this Poisson
pencil. This is the reason why  the quantum sphere (hyperboloid) can
be realized in the both ways: its classical counterpart is a symmetric
orbit!

Note that all these structures are quasiclassical counterparts of
$\Uq$-covariant algebras. However, the family of Poisson structures
possessing this property is much bigger. We refer the reader to
\cite{DGS} where such  Poisson brackets are classified. (In fact on
the sphere the bracket $\{\,\,,\,\,\}_{RE}$ becomes the so-called
R-matrix bracket classified earlier in \cite{GP}.)

Let us complete this section with a discussion of the representation
theory of the RE algebra in the quasiclassical case. Since this
algebra is $\Uq$-covariant it is natural to look for representations
of the RE algebra which would be $\Uq$-morphisms.  Let a linear space
$U$ be a $\Uq$-module. Then the space $\End(U)$ can be equipped with
the structure of a $\Uq$-module as well. Moreover, the product in this
algebra is $\Uq$-covariant. We say that a map
$$
\rho: \Lq\to \End(U)
$$
is a representation if it is true in the category of associative
algebras and if $\rho$ is a $\Uq$-morphism.

Note that this definition looks like the definition of representations
in a super-category.

We say that a representation $\rho$ of the algebra $\Lq$ is {\em
generic} if
$$
\rho(\Tr_q(L))=a\,\Id,\quad a\not=0.
$$
On the contrary, if $a=0$ we say that the representation is {\em
exceptional}.

We can say nothing about exceptional representations. However, the
theory of generic finite-dimensional representations can be
constructed as follows.

Consider the map
$$
\gah:\Lq \to k,\,\, \gah(l_i^j)=h \de_i^j,\,\, h\in k\;.
$$
It is a $U_q(sl(n))$-morphism. Let us replace $L$ in the defining
relations of the RE algebra by $L-\gah(L)$ (i.e., $l_i^j$ by
$l_i^j-\gah(l_i^j)$). Then these relations turn into
\beq
RL_1RL_1-L_1RL_1R=\h(RL_1-L_1R),\,\, \h=h(q-q^{-1}). \label{REh}
\eeq
The algebra defined by \r{REh} is just the algebra $\Uhq$
mentioned above.

By imposing the condition $\Tr_q(L)=0$ we get a q-analogue of the
algebra $\Ush$. Denote this two parameter algebra $\Ushq$. Since this
algebra also is $\Uq$-covariant we consider its representations in the
same sense as above.

It is evident from the construction that the  representations of the
RE algebra such that
$$
\rho(\Tr_q(L))=-\gah(\Tr_q(L))\,\Id
$$
and those of the algebra $\Ushq$ are in the one-to-one correspondence.
(Let us note that $\gah(\Tr_q(L))\not=0$ for a generic $q$.)

The representations of the latter algebra can be constructed in the
quasiclassical case via the embedding
\beq
\Ushq\hookrightarrow\Uq
\label{emb}
\eeq
realized in \cite{LS}. More precisely, embedding \r{emb} differs from
that of \cite{LS} by a factor which is the Casimir element but it
becomes a number if the QG $\Uq$ is represented in an irreducible
module. So, by rescaling the generators of the algebra $\Ushq$ we can
get embedding \r{emb} represented in a module. Thus, any
irreducible representation of the QG $\Uq$ gives rise to that of
$\Ushq$.

Note that in a non-quasiclassical case the above method is not valid.
Another way to develop representation theory of (non-commutative)
quantum hyperboloid was suggested in \cite{DGR}. That method is valid
in all cases mentioned in remark 2. It is based on the fact that if
$\End(U)$ is multiplicity free the map $V\to{\rm End}\,(U)$ being a
morphism is unique up to a factor. If $U$ is the fundamental module
($U={\bf V}$ in the terms of remark 2) $V$ can be identified with the
traceless component of $\End(U)$.

This observation allows us to consider  the problem of defining an
involution operator in the algebra $\ahq$ from a new viewpoint. The
traditional approach begins with definition of this operator on the
algebra $k(G_q)$ (and on its dual object) by means of some
compatibility condition with the coproduct. However, in our
realization of the quantum sphere any Hopf structures (habitual or
braided) are irrelevant. From the point of view of above realization
of the quantum varieties, the problem can be reformulated as follows:
what are desired properties of an involution operator in the space
$\End(V)$, $V$ being an object of a tensor (or quasitensor) category?
(Note that in a super-category the classical condition $(ab)^*=b^*a^*$
is replaced by its super-analogue.) A case of an involutive $R$ was
considered in  \cite{GRZ}. In \cite{DGR} some compatibility condition
with a q-analogue of Lie bracket  was suggested.

\section{Quantum bundles on quantum sphere}

In this section we consider a particular case of the RE algebra
related to the QG $\uq$.
Thus, we get a new realization of the quantum sphere.

In this case the algebra $\Lq$ is generated by four elements
$l_1^1,\,l_1^2,l_2^1,\,l_2^2$.
The space $\span(l_i^j)$ is a direct sum of the trivial component
generated by
$\Tr_q(L)
$ and $V=\Ker \gah$.
By imposing the relation $\Tr_q(L)= 0$ on the RE algebra we get an
algebra generated by the space $V$. Introducing one more constraint
$\det_q(L)=c\not=0$ where $\det_q(L)$ is the quantum determinant we
get just another realization of the quantum sphere (hyperboloid).
Therefore the above numerical polynomial $\oP$ becomes
\beq
\oP=L^2+c_2 \Id
\label{from}
\eeq
with an appropriate non-trivial factor $c_2$.

Thus, we have realized the quantum sphere as a quotient of the RE
algebra. In a similar way the non-commutative analogue $\ahq$ of the
quantum sphere can be realized as a quotient of the algebra
$U(sl(n)_{\h,q})$. The corresponding polynomial $\oP$ has the same
form \r{from} (cf. \cite{GS}).

This realization of the quantum sphere is very useful for definition
of line bundles on it. Let us explain this.

First, we will construct the line bundles related to the fundamental
$\uq$-module $V=V_{1/2}$ (we call these line bundles {\em basic}). Let
$$
M=V\ot \ahq
$$
be a free right $\ahq$-module. It is a $\uq$-module too. Consider its
submodule $\Mn$ generated by the coordinates of the vector
\beq
v\triangleleft L-\nu v,\quad v=(v_1,\, v_2)\in V\;.
\label{act}
\eeq
By $v\triangleleft L$ we mean the vector whose the $j$-th coordinate
is $v_i\ot l^i_j, \, 1\leq i,j\leq 2$.

\begin{proposition} \cite{GS} The quotient-module $\oM_{\nu}=M/\Mn$ is
not trivial iff $\nu$ is a root of the equation $\oP(\nu)=0$. Let
$\nu_1$ and $\nu_2$ be the roots. Then the operators
$$
P_1={(L-\nu_1\Id)\over
(\nu_2-\nu_1)},\qquad P_2={(L-\nu_2\Id)\over(\nu_1-\nu_2)}
$$
are projectors. By definition these operators act on the $\ahq$-module
$M=V\ot\ahq$ as follows
$$
(v_i\ot f^i)\triangleleft L=(v\triangleleft L)_i\ot f^i,\qquad
\forall\;f^i\in
\ahq\;.
$$
Moreover, $P_1+P_2=\Id$ and the quotient $\oM_{\nu_i}$ can be
identified with $\Im\,P_i,\,\, i=1,\,2$. That is the modules
$\oM_{\nu_i},\,\, i=1,2$ are projective.
\end{proposition}

If $q=1$ these modules correspond in the framework of the Serre-Swan
scheme to the "basic line bundles" on the sphere. We want to stress
that in comparison with a construction of these quantum line bundles
in \cite{BM}, \cite{HM} our approach  does not use of any description
of the
quantum sphere via the QG $k(SL_q(2))$ (or what is the same
$k(SU_q(2))$.
Moreover, our construction has an evident generalization to orbits
related to other QG $\Uq, \, n>2$ and is valid in non-quasiclassical
cases as well (cf. \cite{GS}).
The passage to other modules (called  {\em derived}) can be  realized
by an extension of the matrix $L$ to other $\uq$-modules (hereafter
only the algebra $k(S^2_q)$ is considered).

The problem splits into two parts. The first question is: what is a
reasonable way of constructing such an extension?  In the classical
case the extension to the module $V^{\ot k}\ot k(S^2_q)$ can be
introduced by means of  the Leibniz rule
$$
(u\ot v)\triangleleft L=u\ot(v \triangleleft L)+(u \triangleleft L)\ot
v\qquad u, v \in V
$$
if $k=2$ and similarly for $k>2$. The second part of the problem is to
find the CH identity for such an extension of the matrix $L$. In the
classical case the corresponding numerical polynomials can be found
directly from $\oP$ (cf. \cite{GS}).

However, in the quantum case we meet the following difficulty. If we
extended the matrix $L$ by 
the Leibniz rule as above we would be
unable to find the CH identity for the extended matrix (unless the
Hecke symmetry is involutary). An explanation of this fact is
presented in \cite{GS}. 

In this paper we suggest another way to construct an extension of the
matrix $L$.
Define the matrix
\beq
L_+^{(k)}=P_+^{(k)}L_1 P_+^{(k)} \label{ext}
\eeq
where
$$
P_+^{(k)}: V^{\ot k}\to\Sym^k(V)
$$
is the full $q$-symmetrizer of the $k$-th tensor power of $V$. Its
explicit form can be found in \cite{GPS}. 

The action of thus
introduced matrix on the module $V^{\ot k}\ot k(S^2_q)$ is defined as
follows. For any vector 
$$v_{i_1}\otimes\dots\otimes v_{i_k}\otimes
f_j$$
 of this module, the first projector $P_+^{(k)}$ in \r{ext}
extracts its $q$-symmetric component. Then the matrix $L_1$ acts on
the {\em first} element only of each term of this symmetric component
in accordance with the definition of the above Proposition 1.
And, at last, the action of the second projector $P_+^{(k)}$ in
\r{ext} restores the full $q$-symmetry. That is we have an action
$$L_+^{(k)}:\, v^{\otimes k}\otimes k(S^2_q)\rightarrow \Sym^k(V)\ot
k(S^2_q).$$ 

This method is motivated by the fact that in the classical
case such an extension of the matrix $L$  and that arising from the
Leibniz rule are equivalent (up to a non-essential factor) upon
restriction to the symmetric component. Note that for the
construction of line bundles only symmetric component is
relevant.
This way to extend the matrix $L$ enables us to find the CH identity
for the matrix $L_+^{(k)}$ if $k=2$ in the case related to the quantum
sphere (if $k>2$ the problem is still open). More precisely we have
the following.

\begin{proposition} \cite{GS}
If the CH identity for the matrix $L$ is of the form
$$
L^2-aL+b\,\id =0\qquad a=\nu_1+\nu_2,\ b=\nu_1\nu_2\;.
$$
then the matrix $L_+$  defined by \r{ext} $(k=2)$ obeys the CH
identity of the form:
\beq
L_{+}^3-a(1+\frac{q^{- 1}}{2_q})\,L_{+}^2 + (a^2\frac{q^{-1}}{2_q}
- b)\,L_{+}+ab\frac{q^{-1}}{2_q}\,\id =0.
\label{CH-new}
\eeq
\end{proposition}

Remark that if $a=0$ (in particular, this case corresponds  to the
algebra
$k(S^2_q)$) relation \r{CH-new} becomes
$$
L_+^3-b\,L_+=0
$$
and it does not depend on $q$.

\end{document}